\magnification=1200

\catcode`\@=11 
 
\def\nolabels{\def\wrlabel##1{}\def\eqlabel##1{}\def\reflabel##1{}}
\def\writelabels{\def\wrlabel##1{\leavevmode\vadjust{\rlap{\smash%
{\line{{\escapechar=` \hfill\rlap{\sevenrm\hskip.03in\string##1}}}}}}}%
\def\eqlabel##1{{\escapechar-1\rlap{\sevenrm\hskip.05in\string##1}}}%
\def\thlabel##1{{\escapechar-1\rlap{\sevenrm\hskip.05in\string##1}}}%
\def\reflabel##1{\noexpand\llap{\noexpand\sevenrm\string\string\string##1}}}
\nolabels
\global\newcount\secno \global\secno=0
\global\newcount\meqno \global\meqno=1
\global\newcount\mthno \global\mthno=1
\global\newcount\mexno \global\mexno=1
\global\newcount\mquno \global\mquno=1
\global\newcount\tblno \global\tblno=1
\def\newsec#1{\global\advance\secno by1 
\global\subsecno=0\xdef\secsym{\the\secno.}\global\meqno=1\global\mthno=1
\global\mexno=1\global\mquno=1\global\figno=1\global\tblno=1

\bigbreak\medskip\noindent{\bf\the\secno. #1}\writetoca{{\secsym} {#1}}
\par\nobreak\medskip\nobreak}
\xdef\secsym{}
\global\newcount\subsecno \global\subsecno=0
\def\subsec#1{\global\advance\subsecno by1 \global\subsubsecno=0
\xdef\subsecsym{\the\subsecno.}
\bigbreak\noindent{\bf\secsym\the\subsecno. #1}\writetoca{\string\quad
{\secsym\the\subsecno.} {#1}}\par\nobreak\medskip\nobreak}
\xdef\subsecsym{}
\global\newcount\subsubsecno \global\subsubsecno=0
\def\subsubsec#1{\global\advance\subsubsecno by1
\bigbreak\noindent{\it\secsym\the\subsecno.\the\subsubsecno.
                                   #1}\writetoca{\string\quad
{\the\secno.\the\subsecno.\the\subsubsecno.} {#1}}\par\nobreak\medskip\nobreak}
\global\newcount\appsubsecno \global\appsubsecno=0
\def\appsubsec#1{\global\advance\appsubsecno by1 \global\subsubsecno=0
\xdef\appsubsecsym{\the\appsubsecno.}
\bigbreak\noindent{\it\secsym\the\appsubsecno. #1}\writetoca{\string\quad
{\secsym\the\appsubsecno.} {#1}}\par\nobreak\medskip\nobreak}
\xdef\appsubsecsym{}
\def\appendix#1#2{\global\meqno=1\global\mthno=1\global\mexno=1
\global\figno=1\global\tblno=1
\global\subsecno=0\global\subsubsecno=0
\global\appsubsecno=0
\xdef\appname{#1}
\xdef\secsym{\hbox{#1.}}
\bigbreak\bigskip\noindent{\bf Appendix #1. #2}
\writetoca{Appendix {#1.} {#2}}\par\nobreak\medskip\nobreak}
%
%
\def\eqnn#1{\xdef #1{(\secsym\the\meqno)}\writedef{#1\leftbracket#1}%
\global\advance\meqno by1\wrlabel#1}
\def\eqna#1{\xdef #1##1{\hbox{$(\secsym\the\meqno##1)$}}
\writedef{#1\numbersign1\leftbracket#1{\numbersign1}}%
\global\advance\meqno by1\wrlabel{#1$\{\}$}}
\def\eqn#1#2{\xdef #1{(\secsym\the\meqno)}\writedef{#1\leftbracket#1}%
\global\advance\meqno by1$$#2\eqno#1\eqlabel#1$$}
%
%
\def\thm#1{\xdef #1{\secsym\the\mthno}\writedef{#1\leftbracket#1}%
\global\advance\mthno by1\wrlabel#1}
\def\exm#1{\xdef #1{\secsym\the\mexno}\writedef{#1\leftbracket#1}%
\global\advance\mexno by1\wrlabel#1}
%
%
\def\tbl#1{\xdef #1{\secsym\the\tblno}\writedef{#1\leftbracket#1}%
\global\advance\tblno by1\wrlabel#1}
%
\newskip\footskip\footskip14pt plus 1pt minus 1pt 
\def\f@@t{\baselineskip\footskip\bgroup\aftergroup\@foot\let\next}
\setbox\strutbox=\hbox{\vrule height9.5pt depth4.5pt width0pt}
\global\newcount\ftno \global\ftno=0
\def\foot{\global\advance\ftno by1\footnote{$^{\the\ftno}$}}
%
\newwrite\ftfile
\def\footend{\def\foot{\global\advance\ftno by1\chardef\wfile=\ftfile
$^{\the\ftno}$\ifnum\ftno=1\immediate\openout\ftfile=foots.tmp\fi%
\immediate\write\ftfile{\noexpand\smallskip%
\noexpand\item{f\the\ftno:\ }\pctsign}\findarg}%
\def\footatend{\vfill\eject\immediate\closeout\ftfile{\parindent=20pt
\centerline{\bf Footnotes}\nobreak\bigskip\input foots.tmp }}}
\def\footatend{}
%
%
\global\newcount\refno \global\refno=1
\newwrite\rfile
\def\ref{\the\refno\nref}

\def\nref#1{\xdef#1{\the\refno}\writedef{#1\leftbracket#1}%
\ifnum\refno=1\immediate\openout\rfile=refs.tmp\fi
\global\advance\refno by1\chardef\wfile=\rfile\immediate
\write\rfile{\noexpand\item{[#1]\ }\reflabel{#1\hskip.31in}\pctsign}\findarg}
\def\findarg#1#{\begingroup\obeylines\newlinechar=`\^^M\pass@rg}
{\obeylines\gdef\pass@rg#1{\writ@line\relax #1^^M\hbox{}^^M}%
\gdef\writ@line#1^^M{\expandafter\toks0\expandafter{\striprel@x #1}%
\edef\next{\the\toks0}\ifx\next\em@rk\let\next=\endgroup\else\ifx\next\empty%
\else\immediate\write\wfile{\the\toks0}\fi\let\next=\writ@line\fi\next\relax}}
\def\striprel@x#1{} \def\em@rk{\hbox{}}
\def\lref{\begingroup\obeylines\lr@f}
\def\lr@f#1#2{\gdef#1{\ref#1{#2}}\endgroup\unskip}

\def\addref#1{\immediate\write\rfile{\noexpand\item{}#1}} 
\def\startrefs#1{\immediate\openout\rfile=refs.tmp\refno=#1}
\def\xref{\expandafter\xr@f}\def\xr@f[#1]{#1}
\def\refs#1{[\r@fs #1{\hbox{}}]}
\def\r@fs#1{\edef\next{#1}\ifx\next\em@rk\def\next{}\else
\ifx\next#1\xref #1\else#1\fi\let\next=\r@fs\fi\next}
%

%
 \newwrite\ffile\global\newcount\figno \global\figno=1
%
%
\def\fig{\the\figno\nfig}
\def\nfig#1{\xdef#1{\secsym\the\figno}%
\writedef{#1\leftbracket \noexpand~\the\figno}%
\ifnum\figno=1\immediate\openout\ffile=figs.tmp\fi\chardef\wfile=\ffile%
\immediate\write\ffile{\noexpand\medskip\noexpand\item{Figure\ \the\figno. }
\reflabel{#1\hskip.55in}\pctsign}\global\advance\figno by1\findarg}
\def\xfig{\expandafter\xf@g}\def\xf@g \penalty\@M\ {}
\def\figs#1{figs.~\f@gs #1{\hbox{}}}
\def\f@gs#1{\edef\next{#1}\ifx\next\em@rk\def\next{}\else
\ifx\next#1\xfig #1\else#1\fi\let\next=\f@gs\fi\next}
%
%
\newwrite\lfile

{\escapechar-1\xdef\pctsign{\string\%}\xdef\leftbracket{\string\{}
\xdef\rightbracket{\string\}}\xdef\numbersign{\string\#}}

\def\writestop{\def\writestoppt{\immediate\write\lfile{\string\pageno%
\the\pageno\string\startrefs\leftbracket\the\refno\rightbracket%
\string\def\string\secsym\leftbracket\secsym\rightbracket%
\string\secno\the\secno\string\meqno\the\meqno}\immediate\closeout\lfile}}
\def\writestoppt{}\def\writedef#1{}

\def\seclab#1{\xdef #1{\the\secno}\writedef{#1\leftbracket#1}\wrlabel{#1=#1}}

\def\subseclab#1{\xdef #1{\secsym\the\subsecno}%
\writedef{#1\leftbracket#1}\wrlabel{#1=#1}}
\def\appsubseclab#1{\xdef #1{\secsym\the\appsubsecno}%
\writedef{#1\leftbracket#1}\wrlabel{#1=#1}}
\def\subsubseclab#1{\xdef #1{\secsym\the\subsecno.\the\subsubsecno}%
\writedef{#1\leftbracket#1}\wrlabel{#1=#1}}
\newwrite\tfile \def\writetoca#1{}
\def\leaderfill{\leaders\hbox to 1em{\hss.\hss}\hfill}
\def\writetoc{\immediate\openout\tfile=toc.tmp
   \def\writetoca##1{{\edef\next{\write\tfile{\noindent ##1
   \string\leaderfill {\noexpand\number\pageno} \par}}\next}}}

%
\catcode`\@=12 
%
%
%
%
%
\def\dbend{{\manual\char127}}
\def\d@nger{\medbreak\begingroup\clubpenalty=10000
    \def\par{\endgraf\endgroup\medbreak} \noindent\hang\hangafter=-2
    \hbox to0pt{\hskip-\hangindent\dbend\hfill}\ninepoint}
\outer\def\danger{\d@nger}

\def\darr#1{\raise1.5ex\hbox{$\leftrightarrow$}\mkern-16.5mu #1}

%
%
\def\al{\alpha}

  \def\Ga{\Gamma}
\def\de{\delta}

\def\th{\theta}

\def\la{\lambda} \def\La{\Lambda}
\def\rh{\rho}

\def\ph{\phi}    
\def\ch{\chi}

%
%

%

%
%

\def\cO{{\cal O}}

\def\cW{{\cal W}}

\def\proof{\noindent {\it Proof:}\ }
\def\Box{\hbox{$\rlap{$\sqcup$}\sqcap$}}

%

%
%
\def\amsyes{y }

\def\answ{y }

\ifx\answ\amsyes
\input amssym.def


\def\CC{{\Bbb C}}
\def\ZZ{{\Bbb Z}}
\def\NN{{\Bbb N}}

\def\bfg{{\frak g}}  \def\hsltw{{\widehat{\frak{sl}}_2}}
\def\bfh{{\frak h}}

\def\hg{{\widehat{\frak g}}}
\def\whg{\hg}
\else
\def\ZZ{{Z\!\!\!Z}}              
\def\CC{{I\!\!\!\!C}}
\def\NN{{I\!\!N}}

\def\bfg{{\bf g}}
\def\bfh{{\bf h}}

\def\hg{\hat{\bf g}}

  \def\hsltw{\widehat{s\ell_2}}

\fi
%
%
\def\AdM#1{Adv.\ Math.\ {\bf #1}}

\def\CMP#1{Comm.\ Math.\ Phys.\ {\bf #1}}

\def\NPB#1{Nucl.\ Phys.\ {\bf B#1}}
\def\PLB#1{Phys.\ Lett.\ {\bf {#1}B}}

\def\PRD#1{Phys.\ Rev.\ {\bf D#1}}

%

%
%
\def\SMu{\hbox{\lower 3pt\hbox{ \epsffile{su10.eps}}}}
\def\SMs{\hbox{\lower 3pt\hbox{ \epsffile{ss10.eps}}}}
\def\SMd{\hbox{\lower 3pt\hbox{ \epsffile{sd10.eps}}}}

\def\SMS{\leavevmode\vadjust{\rlap{\smash%
{\line{{\escapechar=` \hfill\rlap{\hskip.3in%
                 \hbox{\lower 2pt\hbox{\epsffile{sd10.eps}}}}}}}}}}
\def\SMH{\leavevmode\vadjust{\rlap{\smash%
{\line{{\escapechar=` \hfill\rlap{\hskip.3in%
                 \hbox{\lower 2pt\hbox{\epsffile{su10.eps}}}}}}}}}}
%
%
\def\LW#1{\lower .5pt \hbox{$\scriptstyle #1$}}
\def\LWr#1{\lower 1.5pt \hbox{$\scriptstyle #1$}}
\def\LWrr#1{\lower 2pt \hbox{$\scriptstyle #1$}}
\def\RSr#1{\raise 1pt \hbox{$\scriptstyle #1$}}

%

\hfuzz=10pt
\nopagenumbers
\pageno=0
%
\def\ch{{\rm ch}}
\def\chv{{\rm ch}^V}
\def\La{\Lambda}

%
%
\line{}
\vskip3cm
\centerline{{\bf COSET CONSTRUCTION FOR WINDING SUBALGEBRAS }}
\centerline{{\bf AND APPLICATIONS}}
\vskip1cm

\centerline{Peter BOUWKNEGT\footnote{$^\ddagger$}{Supported by an
Australian Research Council (QEII) Fellowship}}
\bigskip

\centerline{\sl Department of Physics and Mathematical Physics}
\centerline{\sl University of Adelaide}
\centerline{\sl Adelaide, SA~5005, AUSTRALIA}
\medskip
\vskip1.5cm

\centerline{\bf ABSTRACT}\medskip
{\rightskip=1cm 
\leftskip=1cm 
\noindent 
In this paper we review the coset construction for 
winding subalgebras of affine Lie
algebras.  We classify all cosets of central charge $\widehat c<1$ 
and calculate their branching rules.  The corresponding 
character identities give certain `doubling formulae' for the 
affine characters.  We discuss some applications
of our construction, in particular we find a simple proof of a
crucial identity needed for the computation of the level-$2$ 
root multiplicities of the hyperbolic Kac-Moody algebra $E_{10}$.}
\vskip2truecm


\vfil
\line{ADP-96-39/M48 \hfil}
\line{{{\tt q-alg/9610013}}\hfil}

\eject


\lref\BBSS{
F.A.~Bais, P.~Bouwknegt, K.~Schoutens and M.~Surridge,
{\it Extensions of the Virasoro algebra constructed from 
Kac-Moody algebras using higher order Casimir invariants},
\NPB{304} (1988) 348-370; {\it ibid.},
{\it Coset construction for extended Virasoro algebras},
\NPB{304} (1988) 371-391.}

\lref\BH{
K.~Bardak\c{c}i and M.B.~Halpern, {\it New dual quark models},
\PRD{3} (1971) 2493; M.B.~Halpern, {\it The two faces of a dual 
pion-quark model}, \PRD{4} (1971) 2398.}

\lref\Bou{
P.~Bouwknegt, {\it On the construction of modular invariant
partition functions}, \NPB{290} (1987) 507-526.}

\lref\BN{ 
P.~Bouwknegt and W.~Nahm, {\it Realizations of the 
exceptional modular invariant $A^{(1)}_1$ partition functions},
\PLB{184} (1987) 359-362.}

\lref\BG{ 
P.~Bowcock and P.~Goddard, {\it Coset constructions and extended conformal
algebras}, \NPB{305} (1988) 685-709.}

\lref\FL{
V.A.~Fateev and S.L.~Lukyanov, {\it Additional symmetries and 
exactly-soluble models in two-dimensional conformal field theory},
Sov.\ Sci.\ Rev.\ A.\ Phys.\ {\bf 15} (1990) 1.}

\lref\GN{
R.W.~Gebert and H.~Nicolai, 
{\it On $E_{10}$ and the DDF construction}, 
\CMP{172} (1995) 571-622, ({\tt hep-th/9406175}); {\it ibid.},
{\it $E_{10}$ for beginners}, 
in the proceedings of
``G\"ursey Memorial Conference I: Strings and Symmetries,'' {\it eds.}
G.~Akta\c{s} et al., Istanbul, June '94, pp.197-210, 
(Springer Verlag, Berlin, 1995), ({\tt hep-th/9411188}).}

\lref\GKO{
P.~Goddard, A.~Kent and D.~Olive, {\it Virasoro algebras and 
coset space models}, \PLB{152} (1985) 88.}

\lref\Kac{ 
V.G.~Kac, {\it Infinite dimensional Lie algebras}, 3rd edition,
(Cambridge University Press, 1990).} 

\lref\KMW{ 
V.G.~Kac, R.V.~Moody and M.~Wakimoto, {\it On $E_{10}$},
in ``Differential geometrical methods in theoretical physics,'' proc.\
NATO advanced research workshop, 16th international conference, Como,
eds.\ K.~Bleuler and M.~Werner, pp.109-128 (Kluwer, Dortrecht, 1988).}

\lref\KW{ 
V.G.~Kac and M.~Wakimoto, {\it Modular and 
conformal invariance constraints in representation theory of affine 
algebras}, \AdM{70} (1988) 156-236.}

\lref\KWa{
V.G.~Kac and M.~Wakimoto, {\it Branching functions for winding subalgebras 
and tensor products}, Acta Appl.\ Math.\ {\bf 21} (1990) 3.}

\lref\Nak{
H. Nakajima, {\it Instantons on ALE spaces, quiver varieties, and Kac-Moody
algebras}, Duke Math.\ J.\ {\bf 76} (1994) 365-416.}

\lref\Vaf{
C.~Vafa, {\it Instantons on D-branes}, \NPB{463} (1996) 435-442,
({\tt hep-th/9512078}).}

\baselineskip=13pt

\footline{\hss \tenrm -- \folio\ -- \hss}

Let $\bfg$ be a simple, finite-dimensional Lie algebra.  Consider 
the (untwisted) affine Lie algebra $\hg_k$ defined by
$\hg_k = ( \bfg \otimes \CC[t,t^{-1}] ) \oplus
\CC k \oplus \CC d$ with commutators
\eqn\eqaa{ \eqalign{
[x(m), y(n)] & ~=~ [x,y](m+n) + k\,m\, (x|y)\, \delta_{m,-n}\,, \cr
[d   , x(n)] & ~=~ -n\, x(n) \,,\cr
[k   , x(n)] & ~=~ [k,d] ~=~ 0 \,,\cr}
}
where we have written $x(n) = x \otimes t^n$.  We follow the conventions 
of [\Kac].  In particular $(\ | \ )$ denotes the Killing form on 
both $\bfg$ and $\bfh^*$, normalized such that $(\th|\th)=2$ for
a long root $\th$ of $\bfg$.  The integrable highest weight modules 
$L(\La)$ of $\hg$ at level $k$ are parametrized by dominant integral
weights $\La$ such that $\sum a_i^\vee (\La,\al_i^\vee) = k$.  

The proof of the following theorem is standard 
\thm\thaa
\proclaim Theorem \thaa.
\item{i.}
For every $j\in \NN$ we have an embedding $\hg_{jk}
\hookrightarrow \hg_{k}$ defined by $x(n) \mapsto x(jn)$.
\item{ii.} Let $L(\La)$ be an integrable highest weight module 
of $\hg_k$ with highest weight vector $v_\La$.  Considered as a $\hg_{jk}$ 
module, $L(\La)$ is integrable (and hence fully reducible).  Moreover, 
$U(\hg_{jk}) \cdot v_\La\cap L(\La)$ is an irreducible integrable $\hg_{jk}$
module.\par

\noindent (Note: $\hg_{jk}$ is called a winding subalgebra 
of $\hg_{k}$ [\KWa].)

As an example, consider $\hg = \hsltw$.  We have, under $\hg_k \hookrightarrow
\hg_1$ ($k\geq2$), e.g., 
\eqn\eqba{
U(\whg_k)\cdot v_{\La_0}\cap  L(\La_0) ~\cong~ L(k\La_0)\,,
}
and
\eqn\eqbb{
U(\whg_k)\cdot v_{\La_1}\cap  L(\La_1) ~\cong~ L((k-1)\La_0+\La_1)\,.
}
Other irreducible modules can be projected out by considering other
maximal vectors (vectors in $L(\La)$ on the Weyl orbit of the highest
weight vector), e.g.
\eqn\eqbc{
U(\whg_k)\cdot v_{r_0\La_0} \cap L(\La_0) ~\cong~ L((k-2)\La_0+2\La_1)\,.
}

We can actually do more than projecting out the irreducible modules.
The $\hg_{k}$ modules decompose with respect to a direct sum 
of $\hg_{jk}$ and a coset Virasoro algebra.  
The construction is a slight modification
of the standard coset construction [\BH,\GKO].  

Recall that the Virasoro algebra, $V_c$, is generated by 
$\{L(m)\,|\,m\in\ZZ\}$
and a central charge $c$ with relations
\eqn\eqab{
[L(m),L(n)] ~=~ (m-n)\, L(m+n) + {\textstyle{c\over24}} \, 
  m(m^2-1)\, \de_{m,-n}\,,
}
Any (positive energy) module of $\hg_k$ can be extended to a 
module of the semi-direct sum $V_c \oplus \hg_k$ by means of the 
affine Sugawara construction [\BH]
\eqn\eqac{
L^{\bfg,k}(m) ~=~ {1\over 2(k+ {\rm h}^\vee)} \sum_{n\in\ZZ} 
  :x^a(m-n)x^a(n): \,,
}
where ${\rm h}^\vee$ denotes the dual Coxeter number of $\bfg$,
and $\{x^a,\, a=1,\ldots,{\rm dim\,}\bfg\}$ is an orthonormal 
basis of $\bfg$ (with respect to $(\ |\ )$).
The generators $L^{\bfg,k}(m)$ satisfy a Virasoro algebra \eqab\
with central charge 
\eqn\eqad{
c^{\bfg}(k) ~=~ {k\, {\rm dim\,}\bfg \over k+ {\rm h}^\vee} \,,
}
while
\eqn\eqak{
[L^{\bfg,k}(m),x(n)] ~=~ -n\, x(m+n)\,.
}
We can identify $-d = L^{\bfg,k}(0) - h(\La)$ where 
\eqn\eqal{
h(\La) ~=~ { (\La,\La+2\rh) \over 2(k+ {\rm h}^\vee)}\,.
}
The following coset construction was first discovered by Kac and Wakimoto
[\KWa]
\thm\thab
\proclaim Theorem \thab.  Let $j\in\NN$.
\item{i.} Define 
\eqn\eqbd{
L^{\bfg,k,j}(m) ~=~ {\textstyle{1\over j}} L^{\bfg,k}(jm) 
  + {\textstyle{c^{\bfg}(k)\over 24}} 
 \left( {\textstyle{j - {1\over j}}} \right) \de_{m,0} \,,
}
where $L^{\bfg,k}(m)$ denotes the Sugawara generator \eqac.  Then
$L^{\bfg,k,j}(m)$ satisfies a Virasoro algebra with central charge 
$c = j c^{\bfg}(k)$.
\item{ii.} The coset generator
\eqn\eqbf{
\widehat{L}^{\bfg,k,j}(m) ~=~ L^{\bfg,k,j}(m) - L^{\bfg,jk}(m) \,.
}
satisfies a Virasoro algebra with central charge 
\eqn\eqbg{
\widehat{c}^{\bfg}(k;j) ~=~ j c^{\bfg}(k) - c^{\bfg}(jk) \,,
}
and commutes with $\hg_{jk}$.\par
\thm\thac
\proclaim Corollary \thac.  The irreducible integrable highest
weight modules $L(\La)$ of $\hg_k$ decompose into a direct sum
of unitary irreducible modules of $V_{\hat c} \oplus
\hg_{jk}$ under $V_{\hat c} \oplus
\hg_{jk}\hookrightarrow \hg_k$.\par

Irreducible highest weight modules $L(h,c)$ of the Virasoro algebra
are parametrized by the central charge $c$ and the conformal dimension
$h$, i.e.\ the eigenvalue of $L(0)$ on the highest weight vector.
For $0<c<1$ these modules are unitary provided $c$ is in the series
\eqn\eqag{
c(m) ~=~ 1 - {6\over m(m+1)}\,,
}
and, for given $m$, the conformal dimension is an element of the Kac table 
\eqn\eqah{
h^{(m)}(r,s) ~=~ { ( (m+1)r - ms )^2 - 1 \over 4m(m+1)}\,,\qquad 
1\leq r \leq m-1\,, 1\leq s \leq m \,.
}
The character $\chv_h(q)$ of the irreducible $V_c$-module $L(h,c)$ 
is defined by
\eqn\eqai{
\chv_h(q) ~=~ {\rm Tr\,}_{L(h,c)}\  q^{L(0) - c/24} \,.
}
For $h=h^{(m)}(r,s)$ we will also write $\chv_{(r,s)}(q)$.
Similarly, the character $\ch_{L(\La)}(z;q)$ of the 
$\hg$-module $L(\La)$ is defined by
\eqn\eqaj{
\ch_{L(\La)}(z;q) ~=~ {\rm Tr\,}_{L(\La)}\ 
  z^h\, q^{L^{\bfg,k}(0) - c^{\bfg}(k)/24} \,,
}
where $z^h \equiv \prod z_i^{h_i}$ and $h_i$ runs over a basis of the 
Cartan subalgebra of $\bfg$.
\thm\thad
\proclaim Theorem \thad.  The following is a complete list of all simple Lie
algebras $\bfg$ such that $\widehat{c}^{\bfg}(k;j)<1$ together with
the branching rules of all integrable highest weight modules 
of $\hg_k$ under $V_{\hat c} \oplus
\hg_{jk} \hookrightarrow \hg_k$ (up to those related by automorphisms 
of $\bfg$):\par

\item{$\bullet$}  
$(A_1^{(1)})_2\hookrightarrow (A_1^{(1)})_1$, $\widehat{c} = c(3) ={1\over2}$
$$ \eqalign{
q^{{1\over8}}\  \ch_{L(\La_0)}(z;q)    & ~=~
  \chv_{(1,2)}(q^2) \ch_{L(2\La_0)}(z;q^2) + \chv_{(2,2)}(q^2) 
  \ch_{L(2\La_1)}(z;q^2)  \cr
q^{{1\over8}}\  \ch_{L(\La_1)}(z;q)   & ~=~ 
  ( \chv_{(1,1)}(q^2) + \chv_{(2,1)}(q^2))  \ch_{L(\La_0+\La_1)}(z;q^2)\cr}
$$

\item{$\bullet$}
$(A_2^{(1)})_2\hookrightarrow (A_2^{(1)})_1$, $\widehat{c} = c(5) ={4\over5}$
$$ \eqalign{
q^{{1\over4}}\  \ch_{L(\La_0)}(z;q)  & ~=~
  ( \chv_{(1,2)}(q^2) + \chv_{(4,2)}(q^2) ) 
  \ch_{L(2\La_0)}(z;q^2) \cr
& ~+~  ( \chv_{(2,2)}(q^2) + \chv_{(3,2)}(q^2) ) 
  \ch_{L(\La_1+\La_2)}(z;q^2) \cr
q^{{1\over4}}\  \ch_{L(\La_1)}(z;q)  & ~=~
  ( \chv_{(1,2)}(q^2) + \chv_{(4,2)}(q^2) ) 
  \ch_{L(2\La_2)}(z;q^2)\cr
& ~+~ ( \chv_{(2,2)}(q^2) + \chv_{(3,2)}(q^2) ) 
  \ch_{L(\La_0+\La_1)}(z;q^2) \cr}
$$

\item{$\bullet$}
$(E_8^{(1)})_2\hookrightarrow (E_8^{(1)})_1$, $\widehat{c} = c(3)={1\over2}$
$$ \eqalign{
q\  \ch_{L(\La_0)}(z;q)    & ~=~
  \chv_{(2,1)}(q^2) \ch_{L(2\La_0)}(z;q^2) 
  + \chv_{(2,2)}(q^2)\ch_{L(\La_1)}(z;q^2)\cr
& ~+~ \chv_{(2,3)}(q^2) \ch_{L(\La_7)}(z;q^2) \cr}
$$

\item{$\bullet$}
$(E_7^{(1)})_2\hookrightarrow (E_7^{(1)})_1$, $\widehat{c} = 
c(4) = {7\over10}$
$$ \eqalign{
q^{{7\over8}}\  \ch_{L(\La_0)}(z;q)  & ~=~ 
   \chv_{(2,1)}(q^2) \ch_{L(2\La_0)}(z;q^2) 
  + \chv_{(2,2)}(q^2) \ch_{L(\La_1)}(z;q^2) \cr
&  ~+~ \chv_{(3,2)}(q^2) \ch_{L(\La_5)}(z;q^2)  
  + \chv_{(4,2)}(q^2) \ch_{L(2\La_6)}(z;q^2) \cr
q^{{7\over8}}\  \ch_{L(\La_6)}(z;q)  & ~=~
   (\chv_{(1,1)}(q^2) 
 + \chv_{(1,4)}(q^2)) \ch_{L(\La_7)}(z;q^2)  \cr
& ~+~   ( \chv_{(1,2)}(q^2)
 + \chv_{(1,3)}(q^2) ) \ch_{L(\La_0+\La_6)}(z;q^2) \cr}
$$

\item{$\bullet$}
$(E_6^{(1)})_2\hookrightarrow (E_6^{(1)})_1$, $\widehat{c} = 
c(6) = {6\over7}$
$$ \eqalign{
q^{{3\over4}}\ \ch_{L(\La_0)}(z;q)  & ~=~
   (\chv_{(2,1)}(q^2) + \chv_{(2,6)}(q^2)) 
  \ch_{L(2\La_0)}(z;q^2)  \cr
& ~+~ (\chv_{(2,2)}(q^2)+
   \chv_{(2,5)}(q^2))\ch_{L(\La_6)}(z;q^2) \cr
& ~+~ (\chv_{(2,3)}(q^2)+ \chv_{(2,4)}(q^2))
   \ch_{L(\La_1+\La_5)}(z;q^2) \cr
q^{{3\over4}}\ \ch_{L(\La_1)}(z;q)  & ~=~
   (\chv_{(2,1)}(q^2) + \chv_{(2,6)}(q^2)) \ch_{L(2\La_5)}(z;q^2)\cr 
& ~+~ (\chv_{(2,2)}(q^2)+ \chv_{(2,5)}(q^2))\ch_{L(\La_4)}(z;q^2)\cr
& ~+~ (\chv_{(2,3)}(q^2)+ \chv_{(2,4)}(q^2))
   \ch_{L(\La_0+\La_1)}(z;q^2) \cr}
$$

\item{$\bullet$}
$(G_2^{(1)})_2\hookrightarrow (G_2^{(1)})_1$, $\widehat{c} = 
c(9) = {14\over15}$
$$\eqalign{
q^{{7\over20}}\ \ch_{L(\La_0)}(z;q)   & ~=~ 
   (\chv_{(1,2)}(q^2) + \chv_{(8,2)}(q^2)) 
  \ch_{L(2\La_0)}(z;q^2) \cr
&  ~+~ (\chv_{(2,2)}(q^2)+ \chv_{(7,2)}(q^2))
  \ch_{L(\La_1)}(z;q^2) \cr
& ~+~ (\chv_{(3,2)}(q^2) + 
   \chv_{(6,2)}(q^2)) \ch_{L(2\La_2)}(z;q^2) \cr
& ~+~ (\chv_{(4,2)}(q^2) + \chv_{(5,2)}(q^2))
   \ch_{L(\La_0+\La_2)}(z;q^2) \cr
q^{{7\over20}}\ \ch_{L(\La_2)}(z;q)  & ~=~ 
  (\chv_{(1,4)}(q^2) + \chv_{(8,4)}(q^2)) 
  \ch_{L(2\La_0)}(z;q^2) \cr
&  ~+~ (\chv_{(2,4)}(q^2)+ \chv_{(7,4)}(q^2))
  \ch_{L(\La_1)}(z;q^2) \cr
& ~+~ (\chv_{(3,4)}(q^2)+ 
   \chv_{(6,4)}(q^2)) \ch_{L(2\La_2)}(z;q^2)\cr
& ~+~ (\chv_{(4,4)}(q^2)+ \chv_{(5,4)}(q^2))
   \ch_{L(\La_0+\La_2)}(z;q^2) \cr}
$$

\item{$\bullet$}
$(F_4^{(1)})_2\hookrightarrow (F_4^{(1)})_1$, $\widehat{c} = 
c(10) = {52\over55}$
$$\eqalign{
q^{{13\over20}}\ \ch_{L(\La_0)}(z;q)  & ~=~
   (\chv_{(2,1)}(q^2) + \chv_{(2,10)}(q^2)) 
  \ch_{L(2\La_0)}(z;q^2) \cr
& ~+~ (\chv_{(2,2)}(q^2) + \chv_{(2,9)}(q^2)) 
   \ch_{L(\La_1)}(z;q^2) \cr
& ~+~ (\chv_{(2,3)}(q^2) + \chv_{(2,8)}(q^2)) 
   \ch_{L(2\La_4)}(z;q^2) \cr
& ~+~  (\chv_{(2,4)}(q^2) + \chv_{(2,7)}(q^2)) 
   \ch_{L(\La_3)}(z;q^2) \cr
& ~+~ (\chv_{(2,5)}(q^2) + \chv_{(2,6)}(q^2)) 
   \ch_{L(\La_0+\La_4)}(z;q^2) \cr
q^{{13\over20}}\ \ch_{L(\La_4)}(z;q)  & ~=~ 
   (\chv_{(4,1)}(q^2) + \chv_{(4,10)}(q^2)) 
  \ch_{L(2\La_0)}(z;q^2)\cr
& ~+~  (\chv_{(4,2)}(q^2) + \chv_{(4,9)}(q^2))
   \ch_{L(\La_1)}(z;q^2) \cr
& ~+~  (\chv_{(4,3)}(q^2) + \chv_{(4,8)}(q^2))
   \ch_{L(2\La_4)}(z;q^2) \cr
& ~+~ (\chv_{(4,4)}(q^2) + \chv_{(4,7)}(q^2))
   \ch_{L(\La_3)}(z;q^2) \cr
& ~+~  (\chv_{(4,5)}(q^2) + \chv_{(4,6)}(q^2))
   \ch_{L(\La_0+\La_4)}(z;q^2) \cr}
$$
\bigskip

\noindent (Note: the simply-laced cases are also discussed in [\KWa].)

\proof The completeness of the list can be shown either through a case
by case verification or by observing that the coset charge 
$\widehat{c}^{\bfg}(k;j)$ is the same as that of the standard coset
$$
\hg_{jk}~\hookrightarrow~ \underbrace{\hg_k \oplus \ldots \oplus \hg_k}_j
$$ 
(through the 
diagonal embedding) and then applying the results of [\BG].
Since, for given $\widehat{c}^{\bfg}(k;j)$ and $k$ there are only 
a finite number of unitary representations of 
$V_{\hat c}\oplus \hg_{jk}$, the branching
rules are fixed up to a finite number of (integer) factors.  These factors
can be determined by examining the first few terms in the series 
expansions (in $q$) of the characters.  As a check on our calculations we
verified, using the results of [\KW], that the asymptotics of the 
characters are consistent with the obtained branching rules.\Box\bigskip

Note that, even though the coset central charges are the same as those 
of the standard coset construction applied to the diagonal
cosets $\hg_2 \hookrightarrow
\hg_1 \oplus \hg_1$, the branching rules are quite different.  While 
for the diagonal cosets only those $V_c$-modules occur that extend to
modules of the $\cW[\bfg]$-algebra [\FL,\BBSS], in the present case 
typically the complement of those in the Kac table occur.  We expect 
that those modules are modules of an extension of the Virasoro algebra
along the lines of [\BBSS] as well.  

One important application concerns the computation of the root multiplicities
of hyperbolic Kac-Moody algebras or, more generally, generalized Kac-Moody
algebras [\Kac,\KMW,\GN].  As an example,
consider the above branching rule for $\bfg=E_8$.  We find 
\eqn\eqap{
\chv_{(2,2)} (q^2)\, \ch_{L(\La_1)}(z;q^2) ~=~ \left(
  q\, \ch_{L(\La_0)}(z;q) \right)_{\big| {\rm even} }\,,
}
where the subscript `even' denotes that only the even powers 
of $q$ in the series expansion on the right hand side of the 
equation contribute.  Thus, in particular,
\eqn\eqaq{ \eqalign{
q\, b^{\La_1}_{2\La_0}(q^2) & ~=~ q\, b^{\La_1}_{\La_7}(q^2) 
  ~=~ {\ph(q^2)\over \ph(q^4)}\, b^{\La_0}_{\La_0}(q)_{\big| {\rm odd} }\,,\cr
b^{\La_1}_{\La_1}(q^2) & ~=~ {\ph(q^2)\over \ph(q^4)} \,
  b^{\La_0}_{\La_0}(q)_{\big| {\rm even} }\,,\cr}
}
where $b^\La_\la(q)$ denote the Kac-Peterson string functions [\Kac]
(normalized such that $b^\La_\la(q)=\cO(1)$),
and where we have used
\eqn\eqar{
\chv_{(2,2)} (q) ~=~ q^{{1\over16}}\, {\ph(q^2) \over \ph(q)}\,,\qquad
  \qquad \ph(q) ~=~ \prod_{n\geq1} (1-q^n)\,.
}
Or, equivalently,
\eqn\eqas{
b^{\La_1}_{\La_1}(q^2) + q\, b^{\La_1}_{\La_7}(q^2) ~=~ 
  {\ph(q^2)\over \ph(q^4)}\, b^{\La_0}_{\La_0}(q) ~=~
  {\ph(q^2)\over \ph(q^4) \ph(q)^8}\,,
}
which is one of the crucial identities in the computation of the 
level-2 root multiplicities of the hyperbolic Kac-Moody algebra 
$E_{10}$ in [\KMW].  While in [\KMW] equation \eqas\ required a 
lengthy derivation using the modular properties of the characters, we 
see that here it follows quite easily from a coset construction.
It is conceivable that a proper generalization of this 
coset construction along the lines of [\BBSS] will be useful in 
the computation of higher level root multiplicities.

A second application is the construction of modular invariant 
sesquilinear combinations of either $\hg_{jk}$ or $V_c$ characters 
from known modular invariants of $\hg_k$ along the lines of,
e.g., [\BN,\Bou,\KW].  Actually, due to the nature of \eqbd,
these will be invariants of the principal congruence 
subgroup $\Ga(j)$ rather than of the full $PSL(2,\ZZ)$ (and of
$\Ga(2)$ for the examples of Theorem \thad).  We refrain from 
giving a complete list -- the interested reader can easily work out
examples.

Another potential application, which motivated the present work, is
related to the work of Nakajima [\Nak].  It was suggested in [\Vaf]
that the fact that the cohomology of the moduli space of $U(k)$
instantons on the ALE space of type $A_{N}$ carries an
$(A_{N}^{(1)})_k$-module structure [\Nak] can be understood in the
context of heterotic string duality by restricting the $\hg_1$-action
on a free field Fock space to the oscillators mod $k$.  The above
results give a precise description of the decomposition under the
embedding $\hg_k\hookrightarrow \hg_1$.
\bigskip

\noindent {\bf Acknowledgments}\hfil\break
I would like to thank Jim McCarthy and Paul Montague for useful discussions
and Mark Walton for drawing my attention to Ref.\ [\KWa].


\footatend\immediate\closeout\rfile
\baselineskip=14pt{\bigskip\noindent {\bf References}}%
\bigskip{\frenchspacing%
\parindent=20pt\escapechar=` \input refs.tmp\vfill\eject}\nonfrenchspacing


\vfil\eject\end